\documentclass[11pt, reqno]{amsart}

\usepackage{amssymb}

\evensidemargin\oddsidemargin

\begin{document}
\baselineskip=18pt
\setcounter{page}{1}

\newtheorem{Conj}{Conjecture\!\!}
\newtheorem{Theo}{Theorem\!\!}
\newtheorem{Lemm}{Lemma}
\newtheorem{Rem}{Remark}
\newtheorem{Coro}{Corollary\!\!}
\newtheorem{Propo}{Proposition\!\!}

\renewcommand{\theConj}{}
\renewcommand{\theCoro}{}
\renewcommand{\theTheo}{}
\renewcommand{\thePropo}{}

\def\a{\alpha}
\def\b{\beta}
\def\B{{\bf B}}
\def\C{{\bf C}}
\def\cG{{\mathcal{G}}}
\def\cH{{\mathcal{H}}}
\def\cI{{\mathcal{I}}}
\def\cS{{\mathcal{S}}}
\def\UU{{\mathcal{U}}}
\def\ca{c_{\a}}
\def\ka{\kappa_{\a}}
\def\coa{c_{\a, 0}}
\def\cua{c_{\a, u}}
\def\cL{{\mathcal{L}}}
\def\cM{{\mathcal{M}}}
\def\Ea{E_\a}
\def\eps{{\varepsilon}}
\def\esp{{\mathbb{E}}}
\def\Ga{{\Gamma}}
\def\GG{{\bf \Gamma}}
\def\e{{\rm e}}
\def\i{{\rm i}}
\def\L{{\bf L}}
\def\lbd{\lambda}
\def\lcr{\left[}
\def\lpa{\left(}
\def\lva{\left|}
\def\M{{\bf M}}
\def\NN{{\mathbb{N}}}
\def\pb{{\mathbb{P}}}
\def\QQ{{\mathbb{Q}}}
\def\rl{{\mathbb{R}}}
\def\rpa{\right)}
\def\rcr{\right]}
\def\rva{\right|}
\def\W{{\bf W}}
\def\X{{\bf X}}
\def\XX{{\mathcal X}}
\def\YY{{\mathcal Y}}
\def\U{{\bf U}}
\def\V{{\bf V}_\a}
\def\Un{{\bf 1}}
\def\Z{{\bf Z}}
\def\A{{\bf A}}
\def\AA{{\mathcal A}}
\def\hAA{{\hat \AA}}
\def\hL{{\hat L}}
\def\hT{{\hat T}}

\def\claw{\stackrel{d}{\longrightarrow}}
\def\elaw{\stackrel{d}{=}}
\def\qed{\hfill$\square$}

\title{Density solutions to a class of integro-differential equations}

\author[Wissem Jedidi]{Wissem Jedidi}

\address{Department of Statistics and Operation Research, King Saud University, P. O. Box 2455, Riyadh 11451, Saudi Arabia. {\em Email}: {\tt wissem\_jedidi@yahoo.fr}}

\author[Thomas Simon]{Thomas Simon}

\address{Laboratoire Paul Painlev\'e, Universit\'e Lille 1, Cit\'e Scientifique, F-59655 Villeneuve d'Ascq Cedex. {\em Email}: {\tt simon@math.univ-lille1.fr}}

\author[Min Wang]{Min Wang}

\address{Laboratoire Paul Painlev\'e, Universit\'e Lille 1, Cit\'e Scientifique, F-59655 Villeneuve d'Ascq Cedex. {\em Email}: {\tt min.wang@math.univ-lille1.fr}}

\keywords{Beta distribution; Double Gamma function; Fox function; Generalized Gamma convolution; Integro-differential equation; Kr\"atzel function; Riemann-Liouville operator; Stable distribution.}

\subjclass[2010]{26A33; 33E30; 45E10; 60E07}

\begin{abstract} We consider the integro-differential equation ${\rm I}^{\a}_{0+}f= x^m  f$ on the half-line. We show that there exists a density solution, which is then unique and can be expressed in terms of the Beta distribution, if and only if $m> \a.$ These density solutions extend the class of generalized one-sided stable distributions introduced in \cite{WRS} and more recently investigated in \cite{Pk14}. We study various analytical aspects of these densities, and we solve the open problems about  infinite divisibility formulated in \cite{Pk14}.
\end{abstract}

\maketitle

\section{Introduction and statement of the results}

In this paper, we are concerned with the following integro-differential equation
\begin{equation}
\label{Main}
x^m f(x)\; =\; \frac{1}{\Ga(\a)}\,\int_0^x (x-v)^{\alpha-1}\,f(v)\, dv
\end{equation}
on $(0, \infty)$, with $\a > 0$ and $m\in \rl.$ This equation can be written in a more compact way as
$${\rm I}^{\a}_{0+}f\; =\; x^m  f,$$
where ${\rm I}^{\a}_{0+}$ is the left-sided Riemann-Liouville fractional integral on the half-axis. We refer to the comprehensive monograph \cite{KST} for more details on fractional operators and the corresponding differential equations. We are interested in density solutions to (\ref{Main}), that is we are searching for such $f$ satisfying (\ref{Main}) which are also probability densities on $(0, \infty).$ In this framework, the identities (2.1.31) and (2.1.38) in \cite{KST} imply that the auxiliary function $h = {\rm I}^{\a}_{0+}f$ is a solution to the fractional differential equation
\begin{equation}
\label{Main2}
{\rm D}^{\a}_{0+} h\; =\; x^{-m} h,
\end{equation}
where ${\rm D}^{\a}_{0+}$ is the left-sided Riemann-Liouville fractional derivative. This latter equation can be solved in the case $m=1$ in terms of the classical Wright function - see Theorem 5.10 in \cite{KST}, and we will briefly come back to this example in Section 3.

Observe that density solutions to (\ref{Main}) may not exist. If $\a =m$ for example, then (\ref{Main}) becomes
$$f(x)\; =\; \frac{1}{\Ga(\a)}\,\int_0^1 (1-v)^{\alpha-1}\,f(xv)\, dv,$$
and the integral of the right-hand side is infinite if $f$ is non-negative and not identically zero. In this respect, let us also notice that the arbitrary constant $\Gamma(\a)$ in (\ref{Main}) was chosen without loss of generality: if $f_{m, \a}$ is a density solution to (\ref{Main}), then $f_{c, m, \a}(x) = c f_{m, \a}(cx)$ is for every $c > 0$ a density solution to
$$x^m f(x) \;=\; \frac{c^{\a -m}}{\Ga(\a)}\int_0^x (x-v)^{\a-1}f(v) dv.$$
Let us start with a few examples. When $\a = n$ is a positive integer, then (\ref{Main}) becomes an ODE of order $n$ satisfied by the $n-$th cumulative distribution function
$$F_n(x) \; =\; \int_{0< x_1 <\ldots <x_n < x} f(x_1) \, dx_1\ldots dx_n,$$
which is
$$F_n\; =\; x^m F_n^{(n)}.$$

\vspace{2mm}

\begin{itemize}

\item For $\a = 1,$ we solve $F_1 = x^m F_1'$ with $F_1$ bounded and vanishing at zero. This implies that $F_1'$ is a density iff $m > 1$ with $F_1(x) = e^{ -\frac{x^{1-m}}{(m-1)}}$, that is $f_{m,1} = F_1'$ is the density of the Fr\'echet random variable $((m-1)\GG_1)^{\frac{1}{1-m}}$ where, here and throughout, $\GG_t$ denotes a Gamma random variable of parameter $t > 0$, with density
$$f^{}_{\GG_t}(x)\; =\; \frac{x^{t-1} e^{-x}}{\Ga(t)},\qquad x > 0.$$

\item For $\a =2,$ we solve $F_2 = x^m F_2''$ with $F_2$ having linear growth at infinity and vanishing at zero. Supposing $m > 2$ and making the substitution $K(x) = x^{\nu} F_2 ((x/2\nu)^{-2\nu})$ with $\nu = 1/(m-2),$ we obtain Bessel's modified differential equation
$$x^2 K'' \, +\, xK'\, -\, (x^2 + \nu^2) K\; =\; 0,$$
whose solutions satisfying the required properties for $F_2$ are constant multiples of the Macdonald function $K_\nu$. A density solution to (\ref{Main}) is then
$$f_{m,2} (x) \; =\; F_2''(x) \; =\; x^{-m} F_2(x)\; = \; c_\nu\, x^{-3/2 - 1/\nu} K_\nu (2 \nu x^{-1/2\nu}),$$
where $c_\nu$ is the normalizing constant. On the other hand, a computation using e.g. the formula 7.12(23) p.82 in \cite{EMOT} shows that the independent product $\GG_1\times \GG_{\nu +1}$ has density
$$\frac{2\, x^{\frac{\nu}{2}}}{\Ga(\nu+1)}\,K_{\nu}(2 \sqrt{x})\,\Un_{(0,\infty)}(x).$$
By a change of variable, this implies that $f_{m,2}$ is the density of  $((m-2)\sqrt{\GG_1\times \GG_{\nu +1}})^{\frac{2}{2-m}}.$
\end{itemize}

\vspace{4mm}

For $\a \ge 3,$ the resulting ODE's have higher order and do not seem to exhibit any classical special function. In section 3, however, we will see that the density solutions to (\ref{Main}) can be characterized in terms of the Gamma distribution for all integer values of $\a.$

When $\a$ is not a positive integer (\ref{Main}) is a true integro-differential equation, which can be handled via the Laplace transform
$$\cL(\lbd)\; =\; \int_0^\infty e^{-\lbd x} \, f(x)\, dx.$$
In particular, when $m= n$ is a positive integer, the latter satisfies an ODE of order $n$ analogous to the above, which is
$$\cL\; =\; (-1)^n \lbd^\a \cL^{(n)}.$$

\begin{itemize}

\item For $m=1,$ we solve $\cL = - \lbd^\a \cL'$ with $\cL$ a completely monotone (CM) function satisfying $\cL(0) = 1.$ This implies that there is a density solution to (\ref{Main}) iff $\alpha \in (0,1)$ with $\cL(\lbd) = e^{-\frac{\lbd^{1-\a}}{(1-\a)}},$ that is $f_{1,\a}$ is  the density of $(1-\a)^{\frac{1}{\a-1}}\Z_{1-\a}$ where, here and throughout, $\Z_\beta$ is the standard positive $\beta$-stable random variable with Laplace transform
\begin{equation}
\label{LaSta}
\esp[e^{-\lbd \Z_\beta}]\; =\; e^{-\lbd^\beta}, \qquad \lbd \ge 0.
\end{equation}

\item For $m=2,$ we solve $\cL = \lbd^\a \cL''$ with the same restrictions on $\cL.$ Supposing $\a < 2$ and setting $\nu = 1/(2-\a),$ the same reasoning as above leads to
$$\cL(\lbd)\; =\; \frac{2 \nu^\nu}{\Ga(\nu)}\, \sqrt{\lbd} \, K_\nu (2\nu \lbd^{\frac{1}{2\nu}})\; =\; \esp[e^{-\nu^2 \lbd^{\frac{1}{\nu}}\GG_\nu^{-1}}],$$
where the second equality follows again from the formula 7.12(23) p.82 in \cite{EMOT}. For $\a = 1 =\nu,$ we recover the above Fr\'echet density $f_{2,1}(x) = x^{-2} e^{-\frac{1}{x}}.$ For $\a \in (1,2),$ it follows from (\ref{LaSta}) that $f_{2,\a}$ is the density of the independent product $\nu^{2\nu}\,\GG_\nu^{-\nu}\!\times \Z_{\frac{1}{\nu}}.$ For $\a\in (0,1)$, it is not clear from classical  integral formul\ae\, on the Macdonald function that the above $\cL$ is indeed a CM function viz. $f_{2,\a}$ is a density. In Section 3, we will see that $f_{2,\a}$ is for all $\a\in (0,2)$ the density of a certain random variable involving two independent copies of $\Z_{1-\frac{\a}{2}}.$

\end{itemize}

The study of density solutions to (\ref{Main}) for $m$ a positive integer was initiated in \cite{WRS} and then pursued in \cite{Pk14}, where the corresponding random variables are called ``generalized stable". Apart from the classical stable case $m=1,$ these random variables are of interest in the case $\{m = 2, \a\in (0,1)\}$ which is especially investigated in Section 3 of \cite{WRS} and Section 7 of \cite{Pk14}, because of its connections to particle transport along the one-dimensional lattice - see \cite{BSW}. The paper \cite{WRS} takes the point of view of Fox functions and shows that for all $m\in\NN^*, \a\in (0,1)$ there exists a density solution to (\ref{Main})  having a convergent power series representation at infinity and a Fr\'echet-like behaviour at zero - see (2.12) and (2.15) therein. The paper \cite{Pk14} takes the point of view of size-biasing and shows that for all $m\in\NN^*, \a\in (0,m)$ there exists a unique density solution to (\ref{Main}), whose corresponding random variable can be represented in the case $\a\in (m-1,m)$ as a finite independent product involving the random variables $\Z_\beta$ and $\GG_t$ - see Theorems 4.3 and 4.2 therein.

In this paper, we characterize the existence and unicity of density solutions to (\ref{Main}) for all $m \in \rl$ and $\a >0,$ and we obtain a representation of the corresponding random variables as two infinite products involving the Beta random variable $\B_{a,b},$ whose density is recalled to be
$$f^{}_{\B_{a,b}}(x)\; =\; \frac{\Ga(a +b)}{\Ga(a)\Ga(b)}\, x^{a-1}(1-x)^{b-1},\qquad x \in (0,1).$$
Here and throughout, all infinite products are assumed to be independent and a.s. convergent. Our main result reads as follows.

\begin{Theo} The equation (\ref{Main}) has a density solution if and only if $m>\a$. This solution is unique, and it is the density of
$$X_{m, \a}\; \elaw\; \lpa a^{\frac{m-a}{a}} \Ga(\frac{m}{a})\prod\limits_{n=0}^{\infty}\lpa\frac{m+an}{a+an}\rpa\B_{a+an, m-a}\rpa^{-1}\! \elaw\;\,  \lpa\frac{\Gamma(m)}{\Gamma(a)} \prod\limits_{n=0}^{\infty}\lpa\frac{m+n}{a+n}\rpa \B_{1+\frac{n}{a}, \frac{m}{a}-1}\rpa^{-\frac{1}{a}},$$
with the notation $a = m-\a.$
\end{Theo}

The fact that the two above infinite products are actually a.s. convergent is an easy  consequence of the martingale convergence theorem - see the beginning of Section 2.1 in \cite{LS2} and the references therein. The proof of the above theorem relies on the Mellin transform of $f,$ which is by (\ref{Main}) the solution to a functional equation of first order given as (\ref{MGA}) below. This kind of equation has often been encountered in the recent literature, with various point of views - see e.g. \cite{Wb, Pk97, KP, PS3}. If $f$ is assumed to be a density, then (\ref{Main}) or (\ref{MGA}) amount to a random contraction equation
\begin{equation}
\label{RCE}
Y\; \elaw\; \B_{m-\a, \a}\,\times\, {\hat Y}_{m-\a}
\end{equation}
connecting a random variable $Y$ and its size-bias ${\hat Y}_{m-\a},$ with the notations of the beginning of Section 2 in \cite{Pk14}. Our two product representations are then essentially a consequence of Theorem 3.5 in \cite{Pk97} and Lemma 3.2 in \cite{Pk14}. However, these simple representations do not seem to have been observed as yet, see in this respect the bottom of p.208 in \cite{Pk14}.

Throughout, motivated by precise asymptotics analogous to (2.15) in \cite{WRS}, we will also connect the Mellin transform of the solution to (\ref{RCE}) to the double Gamma function $G(z; \tau), z, \tau >0.$ This function, also known as the Barnes function for $\tau =1,$ was introduced in \cite{Bar} as a generalization of the Gamma function. It fulfils the functional equations
\begin{equation}
\label{Bar1}
G(z+1;\tau)\; =\; \Ga(z\tau^{-1})G(z;\tau)\qquad\mbox{and}\qquad G(z+\tau;\tau)\; =\; (2\pi)^{\frac{\tau -1}{2}} \tau^{\frac{1}{2} -z} \Ga(z) G(z;\tau)
\end{equation}
with normalization $G(1;\tau) =1.$ The link between these two functional equations and that of (\ref{MGA}) was thoroughly investigated in \cite{KP} in the framework of L\'evy perpetuities - see Section 3 therein. The normalization also implies
\begin{equation}
\label{Bar2}
G(\tau;\tau)\; =\; \frac{(2\pi)^{\frac{\tau -1}{2}}}{\sqrt{\tau}}
\end{equation}
for every $\tau > 0,$ which will be used henceforth.\\

A consequence of our main result is the solution to open problems related to the  infinite divisibility of the density solutions to (\ref{Main}), recently formulated in \cite{Pk14}. Recall that the law of a positive random variable $X$ is called a generalized Gamma convolution ($X\in\cG$ for short) if there exists a suitably integrable function $a : \rl^+\to\rl^+$ such that
$$X\; \elaw\; \int_0^\infty a(t)\, d\Gamma_t\,,$$
where $\{\Ga_t, \, t\ge 0\}$ is the Gamma subordinator. Equivalently, one has $X\in\cG$ iff its log-Laplace exponent reads
\begin{equation}
\label{GŽgŽ}
-\log \esp[e^{-\lbd X}]\; =\; a\lbd\, +\, \int_0^\infty (1- e^{-\lbd x}) k(x) \frac{dx}{x}\; =\;  a\lbd\, +\, \int_0^\infty \log(1+\lbd u^{-1})  \mu(du)
\end{equation}
with $a \ge 0$ and
$$k(x)\; =\; \int_0^\infty e^{-x u}\, \mu(du)$$
is a CM function, whose Bernstein measure $\mu$ is called the Thorin measure of $X.$ If a random variable $X$ belongs to $\cG,$ then it is self-decomposable ($X\in\cS$ for short), and hence infinitely divisible ($X\in\cI$ for short). An important subclass of $\cG$ is that of hyperbolically completely monotone random variables, which we will denote by $\cH.$ By definition, one has $X\in\cH$ iff $X$ has a positive density $f_X$ on $(0,\infty)$ such that $f_X(uv) f_X(u v^{-1})$ is CM in the variable $v+v^{-1}$ for all $u >0.$  An important property is that
\begin{equation}
\label{IV}
X\,\in\, \cH\;\Leftrightarrow\; X^{-1}\,\in\,\cH,
\end{equation}
showing that $\cH$ is a true subclass of $\cG$ since the inverse of an element of $\cG$ may not even be in $\cI.$ We refer to \cite{Bond} for a classic account on the classes $\cG$ and $\cH$, including all the above facts. See also Chapter VI.5 in \cite{SVH}, and \cite{JRY} for a more recent survey. \\

We now suppose $m > \a$ and denote by $X_{m,\a}$ the positive random variable whose density $f_{m,\a}$ is the unique density solution to (\ref{Main}). The law of $X_{m,\a}$ will be denoted by $\cG(m,\a)$ and called generalized stable with parameters $m$ and $\a,$ whereas the law of $Y_{m,\a} = X_{m,\a}^{-1}$ will be denoted by r$\cG(m,\a),$ in accordance with the terminology of \cite{Pk14}. Observe from (\ref{Main}) that the density $g_{m,\a}(x) = x^{-2} f_{m,\a}(x^{-1})$ of $Y_{m,\a}$ is such that $h_{m,\a} (x) = x^{1-\a} g_{m,\a}(x)$ is a positive solution to
$${\rm I}^{\a}_{-}h\; =\; x^{2\a -m} h,$$
where ${\rm I}^{\a}_{-}$ is the right-sided Riemann-Liouville fractional integral on the half-axis. This dual equation to (\ref{Main}) is the one appearing in a physical context for $m=2$ - see (27) in \cite{BSW}.

 \begin{Coro} With the above notations, one has:

{\em (a)} $X_{m,\a}\in\cG$ for all $m > \a$.

{\em (b)} $X_{m,\a}\in\cH\,\Leftrightarrow \,Y_{m,\a}\in\cI\,\Leftrightarrow\, m \le 2 \a$.

\end{Coro}

Part (a) of the corollary is a generalization of Theorem 5.3 in \cite{Pk14}, solving the open questions formulated thereafter. Besides, by e.g. Theorem III.4.10 in \cite{SVH}, it shows that the density solution $f_{m,\a}$ to (\ref{Main}) is also the unique density solution to
\begin{equation}
\label{Stud}
x f(x)\; =\; \int_0^x k_{m,\a}(x-y)\,f(y) \,dy\,,
\end{equation}
where $k_{m,\a}$ is the CM function associated to $X_{m,\a}$ through (\ref{GŽgŽ}). This latter equation is known as Steutel's integro-differential equation for infinitely divisible densities. Except in the obvious case $m=1$, the link between the two convolution kernels $(x-y)^{\a-1}$ and $k_{m,\a}(x-y)$ is mysterious, and the function $k_{m,\a}$ is not explicit in general. See however Section 3 for an analytical treatment of $k_{m,\a}$ in the cases $m=2\a$ and $m = 2.$ Part (b) gives a characterization of the infinite divisibility of the law r$\cG(m,\a),$ and it is an extension of Theorem 5.1 in \cite{Pk14}. It can also be viewed as a generalization of the main result of \cite{BS3}, which handles the case $m=1.$ As will be observed in Remark \ref{Labr} (a) below, its proof also allows us to solve entirely the open question stated after Theorem 3.2 in \cite{Pk14}.\\

We now turn to the asymptotic behaviour of the densities $f_{m,\a}$ at zero and infinity. This is a basic question for the classical special functions, which is investigated e.g. all along \cite{EMOT}. When $m$ is an integer, the densities $f_{m,\a}$ are Fox functions and in \cite{WRS}, the general results of \cite{Br} are used in order to derive convergent power series representations at infinity, with an exact first order polynomial term, as well as a non-trivial exponentially small behaviour at zero - see (2.20) and (2.21) therein. In general, $f_{m,\a}$ is not a Fox function because its Mellin transform may have poles of infinite order. However, we can show the following estimates, which generalize (2.20) and (2.21) in \cite{WRS}.

\begin{Propo} With the above notation, one has
$$ f_{m,\a}(x)\,\sim\,\frac{x^{\a-m-1}}{\Ga(\a)}\quad \mbox{as $x\to\infty$}\qquad\mbox{and}\qquad f_{m,\a}(x)\,\sim\, c_{m,\a}\, x^{-\frac{m(1+\a)}{2\a}} e^{-(\frac{\a}{m-\a}) x^{\frac{\a-m}{\a}}}\quad\mbox{as $x\to 0,$}$$
with
$$c_{m,\a}\; =\; \frac{(2\pi)^{\frac{m-2}{2}} (m-\a)^{\frac{\a(1-m)}{2(m-\a)}}}{\sqrt{\a}\, G(m,m-\a)}\cdot$$
\end{Propo}

The estimate at infinity is an elementary consequence of (\ref{Main}). The derivation of the estimate at zero, much more delicate, is centered around the exact case $m=2\a$ which corresponds to the Fr\'echet random variable $\GG_\a^{-1}.$ When $m> 2\a,$ the underlying random variable is the exponential functional of a L\'evy process without negative jumps, and we can apply the recent Tauberian results of \cite{PS3}. To handle the case $m< 2\a$ which has fat exponential tails, we perform an induction based on a multiplicative identity in law involving the above $\GG_\a^{-1}.$\\

The next section is devoted to the proof of the three above results. In the last section, we display several remarkable factorizations generalizing those of \cite{Pk14}, and we investigate the corresponding Fox function representations and convergent power series expansions. We also discuss some explicit Thorin measures coming from (\ref{GŽgŽ}) and the behaviour of the laws $\cG(m,\a)$ when $\a \to 0$ and $\a \to m.$

\section{Proofs}
\subsection{Proof of the theorem} We begin with the only if part. Introduce the Mellin transform
$$\mathcal{M}(s) = \int_0^{\infty} x^{-s}f(x)dx$$
which is well-defined for all $s\in\rl$ with possibly infinite values,  since $f$ is non-negative. By Fubini's theorem - see also Lemma 2.15 in \cite{KST}, we readily deduce from (\ref{Main}) the functional equation
\begin{equation}
\label{MGA}
\frac{\mathcal{M}(s+a)}{\mathcal{M}(s)}\; =\; \frac{\Gamma(m+s)}{\Gamma(a+s)}, \qquad s>-a.
\end{equation}
By strict log-convexity of the Gamma function and since $m > m-\a =a,$ we observe that the right-hand side of (\ref{MGA}) is increasing in $s.$ Since $s\mapsto \mathcal{M}(s)$ is also log-convex by H\"older's inequality, the left-hand side of (\ref{MGA}) is non-increasing in $s$ when $a \le 0.$ All of this shows that there is a density solution to (\ref{Main}) only if $a > 0 \Leftrightarrow m > \a.$

We now proceed to the if part. On the one hand, the functional identity (\ref{Bar1}) shows that for all $a\in (0,m)$, the function
\begin{equation}
\label{MMS}
\mathcal{M}(s)\; =\; a^{\frac{(m-a)s}{a}}\times\, \frac{G(m+s,a) G(a,a)}{G(a+s,a)G(m,a)}, \qquad s > -a,
\end{equation}
is a solution to (\ref{MGA}). On the other hand, Proposition 2 in \cite{LS2} implies that this function equals
$$\esp\lcr\lpa a^{\frac{m-a}{a}} \Ga(\frac{m}{a})\prod\limits_{n=0}^{\infty}\lpa\frac{m+an}{a+an}\rpa\B_{a+an, m-a}\rpa^{s}\rcr, \qquad s > -a.$$
This shows that there is a density solution to (\ref{Main}) for $m > \a,$ which is that of $X_{m,\a}$ defined by the first product representation. To obtain uniqueness, we use the same argument as in \cite{Pk14}. If $f$ is a density solution to (\ref{Main}) and if $Y$ is the random variable with density $g(x) = x^{-2} f(x^{-1}),$ we deduce from (\ref{Main}) the identity
$$g(x)\; =\; \frac{\Ga(m)}{\Ga(a)\Ga(m-a)}\int_0^1 t^{a-1} (1-t)^{m-a-1} \lpa \frac{(x t^{-1})^a}{\mathcal{M}-g(a)}g(xt^{-1})\rpa \frac{dt}{t}\,,$$
with the notation
$$\mathcal{M}-g(s) = \int_0^{\infty} x^{-s}f(x)dx\; =\; \esp[Y^s].$$
This translates into the random contraction equation
$$Y\; \elaw\; \B_{a, m-a}\,\times\, {\hat Y}_a\,,$$
where ${\hat Y}_a$ is the size-bias of order $a$ of $Y,$ having density $\frac{x^a g(x)}{\cM-g(a)}.$ By Theorem 3.5 in \cite{Pk97}, the solutions to this random equation are unique up to scale transformation. Since (\ref{MGA}) implies the normalization $\esp[Y^a] = \frac{\Ga(m)}{\Ga(a)},$ we finally obtain the uniqueness of $g,$ and that of $f$ as well.

To conclude the proof, it remains to show the identity in law between the two product representations of $X_{m,\a}.$ This is actually given as Lemma 3.2 in \cite{Pk14}, but we provide here a simple and separate argument. Setting $\cM_a(s) =\cM(as)$ transforms (\ref{MGA}) into
$$\frac{\cM_a(s+1)}{\cM_a(s)}\; =\; \frac{\Gamma(m+as)}{\Gamma(a+as)},$$
whose solution is unique thanks to the main result of \cite{Wb} and the log-convexity of the Gamma function. The functional identity (\ref{Bar1}) shows that this solution is given by
$$\cM_a(s)\; =\; \frac{G(\frac{m}{a}+s,\frac{1}{a}) G(1,\frac{1}{a})}{G(1+s,\frac{1}{a})G(\frac{m}{a},1)}\; =\; \esp\lcr \lpa\frac{\Gamma(m)}{\Gamma(a)} \prod\limits_{n=0}^{\infty}\lpa\frac{m+n}{a+n}\rpa \B_{1+\frac{n}{a}, \frac{m}{a}-1}\rpa^{s}\rcr, \qquad s > -1,$$
where the second equality follows again from Proposition 2 in \cite{LS2}. This completes the proof.

\qed

\subsection{Proof of the Corollary} It is well-known and easy to see from the expression of its density that $\B_{b,c}^{-1}-1\in\cH\subset\cG$ for every $b,c > 0,$ so that $\B_{b,c}^{-1}\in\cG$ as well. The first infinite product representation in the Theorem and the main result of \cite{BJTP} imply that $X_{m,\a}\in\cG,$ which concludes the proof of Part (a).

The first inclusion $X_{m,\a}\in\cH\Rightarrow Y_{m,\a}\in\cI$ of Part (b) is an obvious consequence of (\ref{IV}). As in Theorem 5.1 (b) of \cite{Pk14}, the second inclusion $Y_{m,\a}\in\cI\Rightarrow m \le 2 \a$ follows from well-known bounds on the upper tails of positive ID distributions - see e.g. Theorem III.9.1 in \cite{SVH}, and the small-ball estimate
\begin{equation}
\label{SB}
x^{\frac{\a-m}{\a}}\log \pb[X_{m,\a} < x^{-1}]\; =\; x^{\frac{\a-m}{\a}}\log \pb[Y_{m,\a} > x]\; \to\; -\lpa\frac{\a}{m-\a}\rpa, \qquad x\to \infty.
\end{equation}
When $m$ is a positive integer, the latter estimate is a consequence of (2.15) in \cite{WRS}, taking into account the normalization (2.1) therein. To prove (\ref{SB}) in the general case, we consider the random variable $Z_{m,\a} = (Y_{m,\a})^{\frac{m-\a}{\a}}$ and we study the behaviour of its positive entire moments through the quantities
$$a_n\; =\; \frac{(\esp[(Z_{m,\a})^n])^{\frac{1}{n}}}{n}\; =\; \frac{a}{n} \lpa \frac{G(m+bn,a) G(a,a)}{G(a+bn,a)G(m,a)}\rpa^{\frac{1}{n}},$$
where the second equality follows from (\ref{MMS}), recalling the notation $a = m-\a$ and having set $b =\frac{a}{m-a}\cdot$ By Stirling's formula and the estimate (4.5) in \cite{BK}, we obtain
\begin{equation}
\label{DK}
\lim_{n\to\infty} a_n\; =\;\frac{b}{{\rm e}}
\end{equation}
which, by Lemma 3.2 in \cite{CSY}, implies (\ref{SB}).

In order to show the last inclusion $m\le 2\a \Rightarrow X_{m,\a}\in\cH,$ we will use the argument of the main result in \cite{BS3}. If $m = 2\a,$ then the first product representation and Lemma 3 in \cite{BS3} imply
$$Y_{2\a,\a}\; \elaw\; \a\, \prod\limits_{n=0}^{\infty}\lpa\frac{n+2}{n+1}\rpa\B_{\a(n+1), \a}\; \elaw\; c_\a\, \GG_\a\,,$$
for some normalizing constant $c_\a$ which is here one, since the infinite product has unit expectation. Hence
\begin{equation}
\label{Fresh}
X_{2\a, \a}\,\elaw\, \GG_\a^{-1}\,\in\,\cH.
\end{equation}
If $m < 2\a$ viz. $m > 2a,$ the same argument shows that
\begin{eqnarray*}
Y_{m,\a}& \elaw & a^{\frac{m-a}{a}} \Ga(\frac{m}{a})\prod\limits_{n=0}^{\infty}\lpa\frac{m+an}{a+an}\rpa\B_{a+an, m-a}\\
& \elaw & a^{\frac{m-a}{a}} \Ga(\frac{m}{a})\lpa\prod\limits_{n=0}^{\infty}\lpa\frac{n+2}{n+1}\rpa\B_{a(n+1), a} \rpa \times\lpa \prod\limits_{n=0}^{\infty}\lpa\frac{m+an}{2a+an}\rpa\B_{2a+an, m-2a}\rpa\\
& \elaw & a^{\frac{m-2a}{a}} \Ga(\frac{m}{a})\;\GG_a \; \times\lpa \prod\limits_{n=0}^{\infty}\lpa\frac{m+an}{2a+an}\rpa\B_{2a+an, m-2a}\rpa,
\end{eqnarray*}
which belongs to $\cH$ by Lemma 1 in \cite{BS3}.

\qed

\begin{Rem}
\label{Labr}
{\em (a) The above proof makes it also possible to characterize the infinite divisibility of the class $\cL(a,b,r)$ defined in Section 3 of \cite{Pk14} as the solutions in law to the random contraction equation
$$X\; \elaw\; \B_{a,b}\, \times\, {\hat X}_r,$$
with the notation of Section 2 in \cite{Pk14}. By (3.7) in \cite{Pk14}, these solutions are constant multiples of the infinite product
$$\prod\limits_{n=0}^{\infty}\lpa\frac{a+b+rn}{a+rn}\rpa\B_{a+rn, b},$$
and our argument shows similarly that this product is in $\cH$ as soon as $b\ge r.$ By the second statement of Theorem 3.2 in \cite{Pk14}, this entails the characterization
$$\cL(a,b,r)\,\in\,\cI\;\Leftrightarrow\; \cL(a,b,r)\,\in\,\cH\;\Leftrightarrow\; b \, \ge \, r$$
for every $a > 0,$ providing an answer to the open question stated after Theorem 3.2 in \cite{Pk14}.\\

(b) The second identity in our main result and Example VI.12.21 in \cite{SVH} readily imply that $\log\cG(m,\a)\in\cS$ for every $m > \a.$ Theorem 3.1 in \cite{Pk14} also shows that $\log\cL(a,b,r)\in\cI$ for every $a,b,r >0$, and that $\log\cL(a,b,r)\in\cS$ if and only if the function
$$x\; \mapsto\; \frac{x^a(1-x^b)}{(1-x)(1-x^r)}$$
is non-decreasing on $(0,1),$ a property which neither holds for all $a,b,r >0$ nor seems to be characterized cosily in terms of $a,b,r.$}
\end{Rem}

\subsection{Proof of the Proposition} The asymptotics at infinity is read off immediately in the integro-differential equation (\ref{Main}) itself, since
\begin{eqnarray*}
f_{m,\a}(x) & = & \lpa \int_0^\infty(1-vx^{-1})\Un_{\{v\le x\}} f_{m,\a} (v)\, dv \rpa \frac{x^{\a-m-1}}{\Ga(\a)}\\
& \sim & \frac{x^{\a-m-1}}{\Ga(\a)}\quad \mbox{as $x\to\infty,$}
\end{eqnarray*}
where the estimate follows from dominated convergence and the fact that $f_{m,\a}$ is a density on $(0,\infty).$

\qed

The derivation of the asymptotics at zero is more involved, and we have to consider three cases separately. Observe that at the logarithmic level, the asymptotic was already obtained in (\ref{SB}).

\subsubsection{The case $m = 2\a$} Here, the identity (\ref{Fresh}) implies
$$f_{2\a,\a} (x)\; =\; \frac{x^{-\a-1}}{\Ga(\a)}\, e^{-\frac{1}{x}},$$
which shows the desired asymptotic behaviour in an exact formula, since
$$c_{2\a,\a}\; =\; \frac{(2\pi)^{\a-1} \a^{-\a}}{G(2\a, \a)}\; =\; \frac{(2\pi)^{\frac{\a-1}{2}}}{\sqrt{\a} \Ga(\a) G(\a, \a)}\; =\; \frac{1}{\Ga(\a)}\cdot$$

\qed

\subsubsection{The case $m > 2\a$} We first show the estimate
\begin{equation}
\label{1st}
f_{m,\a}(x)\,\sim\, c\, x^{-\frac{m(1+\a)}{2\a}} e^{-(\frac{\a}{m-\a}) x^{\frac{\a-m}{\a}}}
\end{equation}
for some positive constant $c$ which will be identified afterwards. Recall the notation $a = m-\a$ and introduce the parameter $\beta = \frac{m-a}{a} \in (0,1).$ From (\ref{MMS}) and the first equation in (\ref{Bar1}), we get
$$\phi_{m,\a}(s)\; =\; \frac{\cM(s+1)}{\cM(s)}\; =\; a^\beta\, \frac{\Ga(1+\beta + \frac{s}{a})}{\Ga(1+\frac{s}{a})}\; =\; a^\beta\,\frac{\Ga(1+\beta(1+u))}{\Ga(1+\beta u)}$$
with the notation $s =a\beta u.$ Using e.g. Lemma 1 in \cite{BS1}, this implies
$$\frac{s\,\cM(s+1)}{\cM(s)}\; =\; \psi_{m,\a}(s)$$
where
$$\psi_{m,\a}(s)\; =\; a^\beta\lpa \Ga(\beta +1) s\; +\; \int_{-\infty}^0 (e^{sx} -1 - sx) \lpa \frac{am\beta e^{mx}}{\Ga(1-\beta) (1- e^{ax})^{\beta +2}}\rpa dx\rpa$$
is the Laplace exponent of a L\'evy process without positive jumps $\{L^{(m,\a)}_t, \, t\ge 0\},$ that is $\psi_{m,\a}(s) = \log \esp [e^{s L^{(m,\a)}_1}].$ By the Bertoin-Yor criterion - see Proposition 2 in \cite{BY} and its proof, we deduce
$$X_{m,\a}\; \elaw\; \int_0^\infty e^{-L^{(m,\a)}_t}\, dt.$$
The required estimate will now follow from a recent general result of Patie and Savov on exponential functionals of L\'evy processes without positive jumps. We first write
$$\phi_{m,\a}(s)\; =\; a^\beta\,\Phi_\beta (sa^{-1})\,,$$
where
$$\Phi_\beta(u)\; =\; \lpa\frac{u+\beta}{u}\rpa \frac{\Ga(\beta +u)}{\Ga(u)}\; = \; u^\beta\lpa1 \, +\, \frac{\beta(\beta +1)}{2u}\, +\,  \frac{\beta(\beta^2 -1)(3\beta +2)}{24u^2}\, +\, {\rm O} (u^{-3})\rpa,$$
the expansion being e.g. a consequence of Formul\ae\, (4) and (5') in \cite{ET}. This expansion also shows, after some algebra, that
$$\Phi_\beta^{-1}(u)\; =\;  u^{\frac{1}{\beta}}\, -\, \frac{(\beta +1)}{2}\, +\, {\rm O} (u^{-\frac{1}{\beta}})$$
and, from the concavity of $\Phi_\beta$ and the monotone density theorem, that
$\Phi_\beta'(u) \sim  \beta u^{\beta-1}.$ This implies
$$\phi_{m,\a}^{-1}(s)\; =\; a\,\Phi_\beta^{-1} (sa^{-\beta})\; =\; s^{\frac{1}{\beta}}\, -\, \frac{m}{2}\, +\, {\rm O} (s^{-\frac{1}{\beta}})$$
and
$$(\phi_{m,\a}^{-1})'(s)\; =\;\frac{1}{\phi_{m,\a}'(\phi_{m,\a}^{-1}(s))}\; \sim\; \frac{s^{\frac{1-\beta}{\beta}}}{\beta}\cdot$$
Putting everything together with Formula (5.47) in \cite{PS3}, we finally obtain (\ref{1st}), and it remains to identify the constant $c.$ To do so, we introduce the random variable $U_{m,\a} = \beta X_{m,\a}^{-\frac{1}{\beta}} =\beta Z_{m,\a},$ with density
$$h_{m,\a}(x)\; =\; (\beta x^{-1})^{\beta +1} f_{m,\a} ((\beta x^{-1})^{\beta})\; \sim\; c\,(\beta x^{-1})^{\frac{m(1-\a)}{2(m-\a)}}\, e^{-x}\quad\mbox{as $x\to \infty.$}$$
A standard approximation using Laplace's method and Stirling's formula implies
$$\frac{\esp[U_{m,\a}^n]}{n!}\; \sim\; c\,(\beta n^{-1})^{\frac{m(1-\a)}{2(m-\a)}}\quad\mbox{as $n\to \infty.$}$$
On the other hand, we have
$$\frac{\esp[U_{m,\a}^n]}{n!}\; =\; \frac{\a^n}{n!} \lpa \frac{G(m+bn,a) G(a,a)}{G(a+bn,a)G(m,a)}\rpa\; \sim\; c_{m,\a}\,(\beta n^{-1})^{\frac{m(1-\a)}{2(m-\a)}}\quad\mbox{as $n\to \infty,$}$$
where the estimate follows from (4.5) in \cite{BK}, Stirling's formula, and some algebra. This completes the proof.

\qed

\subsubsection{The case $m < 2\a$} In this case, the small ball estimate (\ref{SB}) shows that $Y_{m,\a}$ does not have exponential moments, so that $X_{m,\a}$ is not distributed as the exponential functional of a L\'evy process without positive jumps, by Proposition 2 in \cite{BY}. Hence, we cannot use the estimate (5.47) in \cite{PS3}. We will first prove (\ref{1st}) via an induction on $n,$ where
\begin{equation}
\label{naa}
(n+1) a\; <\; m\; \le\; (n+2) a.
\end{equation}
The case $n=0$ follows from the previous cases $m\ge 2\a\Leftrightarrow m\le 2a.$ To prove the induction step, we first observe the identity in law
\begin{equation}
\label{FreshFresh}
Y_{m,m-a}\; \elaw\; Y_{m-a, m-2a}\;\times\;\GG_\a,
\end{equation}
which is a consequence of (\ref{MMS}), the second equation in (\ref{Bar1}), and fractional moment identification. The multiplicative convolution formula leads then to
$$f_{m,\a} (x)\; =\; \frac{1}{\Ga(\a)} \int_0^\infty f_{\a, 2\a -m} (xy)\, y^\a \, e^{-y}\, dy.$$
Setting again $b = \frac{a}{m-a} < 1,$ we choose $\delta \in (b, 1)$ and we decompose
\begin{equation}
\label{Landau}
f_{m,\a} (x)\; =\; \frac{1}{\Ga(\a)} \int_0^{x^{-\delta}} f_{\a, 2\a -m} (xy)\, y^\a \, e^{-y}\, dy\; +\; {\rm o} (e^{-x^{-\eta}}),\qquad x\to 0,
\end{equation}
for every $\eta \in (b, \delta),$ where the Landau estimate follows readily from the bounded character of $f_{\a, 2\a -m}.$ To estimate the integral, we use the induction hypothesis on $f_{\a, 2\a -m},$ and the fact that $\delta < 1,$ in order to obtain
\begin{eqnarray*}
\int_0^{x^{-\delta}}\!\!\! f_{\a, 2\a -m} (xy)\, y^\a \, e^{-y}\, dy & \sim & c\, x^{-\frac{\a(1+2\a -m)}{2(2\a-m)}} \int_0^{x^{-\delta}} y^{\frac{\a(2\a -m-1)}{2(2\a-m)}} \, e^{-y - (\frac{2\a-m}{m-\a})(xy)^{\frac{\a-m}{2\a -m}}} \, dy\\
& \sim & c\, x^{-(b +\frac{m+1}{2})} \int_0^{x^{b-\delta}} z^{\frac{\a(2\a -m-1)}{2(2\a-m)}} \, e^{-x^{-b} (z + (\frac{2\a-m}{m-\a})z^{\frac{\a-m}{2\a -m}})} \, dz.
\end{eqnarray*}
Using Laplace's approximation, we deduce that there exists a positive constant ${\tilde c}$ such that
$$\int_0^{x^{-\delta}}\!\!\! f_{\a, 2\a -m} (xy)\, y^\a \, e^{-y}\, dy \; \sim \; {\tilde c}\, x^{-\frac{b + m+1}{2}} \, e^{-b^{-1} x^{-b}}\; =\; {\tilde c}\,x^{-\frac{m(1+\a)}{2\a}} e^{-(\frac{\a}{m-\a}) x^{\frac{\a-m}{\a}}}.$$
By (\ref{Landau}), this completes the proof of (\ref{1st}) by induction. The identification of the constant ${\tilde c}$ is done exactly in the same way as in the case $m > 2\a.$

\qed

\begin{Rem} {\em (a) The derivation of the asymptotics at infinity follows also, in a more complicated way similar to the argument of Theorem 4.4 in \cite{Pk14}, from the behaviour of $\cM(s)$ at its first pole $s=-a.$ More precisely, by (\ref{MGA}), we have
$$\cM(s) \;  \sim \; \lpa \frac{a^{a-m} G(m-a,a)G(a,a)}{G(m,a)}\rpa \times\, \frac{1}{G(a+s,a)}\; =\; \frac{1}{a \Ga(m-a) G(a+s,a)}\qquad\mbox{as $s\downarrow -a,$}$$
where the equality comes from (\ref{Bar2}) and the second equation in (\ref{Bar1}). The latter also imply
$$\frac{1}{G(a+s,a)}\; =\; \frac{(2\pi)^{\frac{a-1}{2}} a^{\frac{1}{2} -a-s} \Ga(s+a)}{G(2a+s,a)}\; \sim\; \frac{a}{a+s}\qquad\mbox{as $s\downarrow -a,$}$$
showing that this first pole is simple and isolated. Putting everything together and using e.g. Theorem 4 in \cite{FGD}, we obtain the required asymptotic
$$f_{m,\a}(x)\; \sim \; \frac{x^{-a-1}}{\Ga(m-a)}\; =\; \frac{x^{\a-m-1}}{\Ga(\a)}\qquad \mbox{as $x\to \infty.$}$$
In principle, the exact expression of $\cM(s)$ and Theorem 4 in \cite{FGD} should make it possible to derive a more complete expansion of $f_{m,\a}$ at infinity. As mentioned in the introduction, an absolutely convergent power series expansion exists when $m$ is an integer, as a consequence of a Fox representation of order $m$ for $g_{m,\a} (x) = x^{-2} f_{m,\a}(x^{-1})$ - see (2.11) and (2.12) in \cite{WRS}. However, since the double Gamma function may have poles of infinite order, the existence of such a convergent power series expansion is delicate in general. We will consider some examples in Section 3.1, revisiting in particular the case when $m$ is an integer.  See also Theorem 3 in \cite{KP} for some results in this vein, which apply to some cases when $m > 2\a$ is not an integer.\\

(b) In the strict stable case $m=1$, our asymptotic at zero reads simply
$$f_{1,\a} (x) \; \sim\; \frac{x^{-\frac{1+\a}{2\a}}}{\sqrt{2\pi\a}}\, e^{-\lpa \frac{\a}{1-\a}\rpa x^{\frac{\a-1}{\a}}}\; =\; \frac{x^{-\frac{2-a}{2(1-a)}}}{\sqrt{2\pi(1-a)}}\, e^{-\lpa \frac{1-a}{a}\rpa x^{\frac{-a}{1-a}}},$$
in accordance with $X_{1,\a}\elaw a^{-\frac{1}{a}} \Z_a$ - see the first order term of (2.4.30) in \cite{IL} and also Theorem 1 in \cite{jedidi}. The latter formula displays actually a complete expansion of the density $f_{1,\a}$ at zero, with non-explicit coefficients. The detailed argument for this expansion, which relies on (\ref{LaSta}), Fourier inversion, and the method of steepest descent, is in the proof of Theorem 2.4.6 in \cite{IL}. In the absence of explicit Laplace transform, a complete expansion at zero for $f_{m,\a}$ seems difficult to derive in general when $m\neq 1.$ \\

(c) The multiplicative identity (\ref{FreshFresh}) has a more general formulation, which is
\begin{equation}
\label{Frech}
Y_{m,\a}\; \elaw\; Y_{q,q-a}\; \times\; Y_{m+a-q,m-q}^{(q-a)}
\end{equation}
for every $q\in(m-a, m),$ with the alternative notation $X^{(t)} = {\hat X}_t.$ Notice that (\ref{Frech}) boils down to (\ref{FreshFresh}) for $q = m-a$ and that, contrary to the self-similar identity
$$X_{1,1-a}\; \elaw\; b^{\frac{a-b}{ab}}\;X_{1,1-b}\; \times\; X_{1,1-\frac{a}{b}}^{\frac{1}{b}}$$
which is valid for every $b\in (a,1),$ it is not a subordination formula. As in Corollary 4 (a) of \cite{LS2}, it can also be shown that $X_{m,\beta}$ is a multiplicative factor of $X_{m,\a}$ for every $0 < \beta < \a < m.$}

\end{Rem}

\section{Further remarks}

\subsection{Some particular factorizations}

In this paragraph we consider three situations where the law $\cG(m,\a)$ has simpler expressions as a finite product involving the Gamma or the positive stable distribution. This expression is derived from rewriting (\ref{MMS}) as a moment of Gamma type, thanks to the concatenation formul\ae\, of (\ref{Bar1}). We refer to \cite{J} for a survey on moments of Gamma type. In our three cases, the density $f_{m,\a}$ is also a Fox $H-$function and we display the convergent power series representations, when it is possible. Throughout, we use again the notations $a = m-\a$ and $X^{(t)} = {\hat X}_t$ in order to have simpler formul\ae. Our reference for  Fox functions is Section 1.12 in \cite{KST}, especially (1.12.1) and (1.12.19) therein.

\subsubsection{The case $\a = n\in\NN$} We have
$$\cM(s)\; =\; a^{\frac{ns}{a}}\;\prod_{i=1}^n \lpa \frac{\Ga(1+ \frac{i-1}{a} + \frac{s}{a})}{\Ga(1+ \frac{i-1}{a})}\rpa, \quad s > -a.$$
This shows that $X_{m,n}$ is a finite independent product of generalized Fr\'echet random variables, as was already observed in the introduction for $\a = 1,2:$ one has
$$X_{m,n}\;\elaw\; \lpa a^n\, \GG^{}_1\,\times\,\cdots\,\times\,\GG_{1 + \frac{n-1}{a}}\rpa^{-\frac{1}{a}}.$$
The Fox function representation of $f_{m, n}$ is then
$$f_{m,n}(x)\; =\; \lpa \frac{a^{\frac{n}{a}}}{ \prod_{i=1}^n \Ga(1+ \frac{i-1}{a})}\rpa H^{0,n}_{n,0}\lcr a^{\frac{n}{a}}x\lva \begin{array}{l} (\frac{-i}{a}, \frac{1}{a})^{}_{i=1, \dots, n}\\ \\ \hline \\ \end{array}\right.\rcr.$$
When $a\not\in\QQ$ or $a\in\QQ$ with $a =\frac{p}{q}$ irreducible and $p\ge n,$ the following convergent power series representation holds:
$$f_{m,n}(x)\; =\; \lpa \frac{a^{\frac{n}{a}+1}}{ \prod_{i=1}^n \Ga(1+ \frac{i-1}{a})}\rpa\sum_{r=1}^n\sum_{k=0}^\infty \frac{(-1)^k a^{-(\frac{rn}{a} + n(k+1))}}{k!}\lpa{\hat\prod}_{j= 1}^n {\textstyle \Ga(\frac{j-r}{a} -k)}\rpa x^{-(r+a(k+1))},$$
where the hat product indicates omission of $j=r.$ For $n=1,$ this simplifies into
$$f_{m,1}(x) \; =\; x^{-a-1}\,\sum_{k\ge 0} \frac{(-1)^k}{k!} \frac{1}{(ax^a)^k}\;=\; x^{-a-1}e^{- \frac{1}{ax^a}}$$
as expected, since $X_{m,1}\elaw (a\GG_1)^{-\frac{1}{a}}.$ For $n=2$ and $\nu = \frac{1}{a}\not\in \NN,$ this simplifies into
\begin{eqnarray*}
f_{m,2}(x) & = & \frac{x^{-a-1}}{\Ga(\nu)}\,\sum_{k\ge 0} \frac{(-1)^k}{k!} \lpa \Ga(\nu-k) +\Ga(-\nu-k) (ax^{a/2})^{-2\nu}\rpa (ax^{a/2})^{-2k}\\
&= & \frac{2 x^{-a-3/2}}{a^{\nu}\Ga(\nu)}\, K_\nu (\frac{2}{ax^{a/2}})
\end{eqnarray*}
as expected, since $X_{m,2}\elaw (a\sqrt{\GG_1\times\GG_{1+\nu}})^{-\frac{2}{a}}$ - see the second example in the introduction. Notice that the representation of $f_{m,2}$ in terms of the Macdonald function $K_\nu$ also holds for $\nu\in\NN,$ but then the convergent series representation has a logarithmic term - see Formula 7.2.5(37) in \cite{EMOT}.

\subsubsection{The case $m=an, n\in\{2,3,\ldots\}$} We have
$$\cM(s)\; =\; \prod_{i=1}^n \lpa \frac{\Ga(ia + s)}{\Ga(ia)}\rpa, \quad s > -a.$$
This shows that $X_{an,a(n-1)}$ is a finite independent product of inverse Gamma random variables, as was already observed in (\ref{Fresh}) for $n = 1:$ one has
$$X_{m,n}\;\elaw\; \lpa \GG^{}_a\,\times\,\cdots\,\times\,\GG_{a(n-1)}\rpa^{-1}.$$
The Fox function representation of $f_{an,a(n-1)}$ is
$$f_{an,a(n-1)}(x)\; =\; \lpa \frac{1}{\prod_{i=1}^{n-1} \Ga(ia)}\rpa H^{0,n-1}_{n-1,0}\lcr x\lva \begin{array}{l} (-ia, 1)^{}_{i=1, \dots, n}\\ \\ \hline \\ \end{array}\right.\rcr.$$
When $n=2$ or $a, \ldots, (n-2)a\not\in\NN,$ the following convergent power series representation holds:
$$f_{an,a(n-1)}(x)\; =\; \lpa \frac{1}{\prod_{i=1}^{n-1} \Ga(ia)}\rpa\sum_{r=1}^{n-1}\sum_{k=0}^\infty \frac{(-1)^k}{k!}\lpa{\hat \prod}_{j= 1}^{n-1} \Ga((j-r)a -k)\rpa x^{-ra-k-1}.$$
For $n=2,$ this simplifies into
$$f_{2a,a}(x)\; =\; \frac{x^{-a-1}}{\Ga(a)}\, \sum_{k\ge 0} \frac{(-1)^k x^{-k}}{k!}\; =\; \frac{x^{-a-1}e^{-\frac{1}{x}}}{\Ga(a)},$$
as expected from (\ref{Fresh}). For $n =3$ and $a\not\in \NN,$ similarly as above we get
$$f_{3a,2a}(x)\; =\;\frac{x^{-a-1}}{\Ga(a)\Ga(2a)}\, \sum_{k\ge 0}\frac{(-1)^k}{k!}\lpa \Ga(a-k) +\Ga(-a-k) x^{-a}\rpa x^{-k}\; =\; \frac{x^{-3a/2-1}}{\Ga(a)\Ga(2a)}\,K_a(2x^{-1/2}),$$
the representation on the right-hand side in terms of the Macdonald function holding for $a\in\NN$ as well.

\subsubsection{The case $m=n$} We have
$$\cM(s)\; =\; a^{\frac{ns}{a}}\;\frac{\Ga(1+\frac{s}{a})}{\Ga(1+s)}\,\times\,\prod_{i=1}^{n-1} \lpa \frac{\Ga(\frac{i+s}{a})}{\Ga(\frac{i}{a})}\rpa\; =\; \lpa \frac{a^{\frac{n}{a}}}{n}\rpa^s\,\times\,\prod_{i=0}^{n-1} \lpa \frac{\Ga(1+\frac{i+s}{a})\Ga(1+\frac{i}{n})}{\Ga(1+\frac{i+s}{n})\Ga(1+\frac{i}{a})}\rpa, \quad s > -a,$$
where the second equality comes from the Legendre-Gauss multiplication formula for the Gamma function. Observe also that the first equality is Theorem 4.1 in \cite{Pk14}, with a different  normalization. This shows that $X_{n,\a}$ is a finite independent product of power transforms of size-biased stable random variables, as was already mentioned in the introduction for $m = 1,2:$ one has
$$X_{n,\a}\;\elaw\; n a^{-\frac{n}{a}}\lpa  \Z^{}_{\frac{a}{n}}\,\times\,\Z^{(\frac{-1}{n})}_{\frac{a}{n}}\,\times\,\cdots\,\times\,\Z_{\frac{a}{n}}^{(\frac{1-n}{n})}\rpa^{\frac{1}{n}}.$$
This factorization may look more satisfactory than that of Theorem 4.2 in \cite{Pk14}, and it is also valid in the full range $\a\in (0,n).$ The Fox function representation of $f_{n, \a}$ is derived similarly as (2.11) in \cite{WRS} - beware again our different normalization: one has
$$f_{n,\a}(x)\; =\; \lpa \frac{a^{\frac{n}{a}-1}}{ \prod_{i=1}^{n-1} \Ga(\frac{i}{a})}\rpa H^{0,n}_{n,1}\lcr a^{\frac{n}{a}}x\lva \begin{array}{l} (1-\frac{i}{a}, \frac{1}{a})^{}_{i=1, \dots, n}\\ \\ (0,1) \\ \end{array}\right.\rcr.$$
When $a\not\in\QQ$ or $a\in\QQ$ with $a =\frac{p}{q}$ irreducible and $p\ge n,$ the following convergent power series representation holds:
$$f_{n, \a}(x)\; =\; \lpa \frac{a^{\frac{n}{a}}}{ \prod_{i=1}^{n-1} \Ga(\frac{i}{a})}\rpa\sum_{r=1}^n\sum_{k=0}^\infty \frac{(-1)^k a^{-(\frac{rn}{a} + nk)}}{k! \Ga(1-r-ak)}\lpa{\hat\prod}_{j= 1}^n {\textstyle \Ga(\frac{j-r}{a} -k)}\rpa x^{-r-ak}.$$
For $n=1,$ this simplifies into
$$f_{1,\a}(x) \; =\; \sum_{k\ge 1} \frac{(-a)^{-k}}{k!\Ga(-ak)}\, x^{-1-ak}$$
as expected from e.g. Theorem 2.4.1 in \cite{IL}, since $X_{1,\a}\elaw a^{-\frac{1}{a}}\Z_a.$ By (2.2.35) in \cite{KST}, the auxiliary function $h_{1,\a} = {\rm I}^{\a}_{0+}f_{1,\a}$ has Laplace transform
$$(\cL h_{1,\a}) (\lbd)\; =\; \lbd^{-\a} \exp \lpa -\frac{\lbd^{1-\a}}{1-\a}\rpa,$$
in accordance with (5.2.143) in \cite{KST}, which leads to (5.2.139) therein, and our above equation (\ref{Main2}) which is for $m=1$ the fractional differential equation (5.2.137) in \cite{KST} with $\lbd = 1$ therein.

In the physically relevant case $n=2$  and for $\nu=\frac{1}{a}\not\in \NN,$ the series representation simplifies into
$$f_{2, \a}(x) \; = \; \frac{1}{x \Ga(\nu)}\,\sum_{k\ge 0} \frac{(-1)^k}{k!\Ga(-ak)} \lpa \Ga(\nu-k) -(1+ak)\Ga(-\nu-k) (ax^{a/2})^{-2\nu}\rpa (ax^{a/2})^{-2k}.$$
Observe the striking formal resemblance with $f_{m,2}$, although no expression in terms of a classical special function seems here available. Notice also that for $\nu\in\NN,$ there is no convergent power series representation for $f_{2,\a}$ in general, save for $\nu=a=\a=1$ where the reduction formula (1.12.43) in \cite{KST} yields
$$f_{2,1}(x)\; =\; H^{0,2}_{2,1}\lcr x\lva \begin{array}{l} (1-i, 1)^{}_{i=1, 2}\\ \\ (0,1) \\ \end{array}\right.\rcr\; =\; H^{0,1}_{1,0}\lcr x\lva \begin{array}{l} (-1, 1)\\ \\ \hline \\ \end{array}\right.\rcr\; =\; x^{-2} e^{-\frac{1}{x}},$$
as again expected from (\ref{Fresh}).

\subsection{Some explicit Thorin measures} As mentioned in the introduction, it follows from Part (a) of the Corollary that the density solutions to (\ref{Main}) are also solution to the Steutel's integro-differential equation (\ref{Stud}), whose convolution kernel $k_{m,\a}(x-y)$ is such that
$$k_{m,\a}(x)\; =\; \int_0^\infty e^{-xu}\mu_{m,\a} (du)$$
is a CM function. In the literature, the measure $\mu_{m,\a}$ is called the Thorin measure associated to the random variable $X_{m,\a}\in\cG$, and we refer to \cite{JRY} - see also Chapter 3 in \cite{Bond} - for more on this topic. From (\ref{GŽgŽ}), the measure $\mu_{m,\a}$ is related to the Laplace transform $\cL_{m,\a}$ of $X_{m,\a}$ via its Stieltjes transform:
$$ \int_0^\infty \frac{\mu_{m,\a}(du)}{u+\lbd}\;=\; -(\log \cL_{m,\a})'(\lbd).$$
Recall that when $m =1,$ we have
$$k_{1,\a}(x)\; =\; \frac{x^{\a-1}}{\Ga(\a)}\; =\; \int_0^\infty e^{-xu}\lpa \frac{\sin(\pi\a)}{\pi u^\a}\rpa du,$$
so that $\mu_{1,\a}$ has a simple explicit density. Let us mention two other cases where $\mu_{1,\a}$ has a more or less explicit density.

\subsubsection{The case $m=2\a$} This case was already discussed at the end of Section 3.2 in \cite{BS1}, but we do it again here for completeness. From (\ref{Fresh}) we have $X_{2\a,\a}\elaw\GG_\a^{-1},$ whose Laplace transform is computed similarly as in the introduction:
$$\cL_{2\a, \a}(\lbd)\; =\; \frac{2\, \lbd^{\frac{\a}{2}}}{\Ga(\a)}\,K_{\a}(2 \sqrt{\lbd}).$$
Using Formul\ae\, 7.11.(25-26) in \cite{EMOT}, we deduce
\begin{equation}
\label{PNT}
-(\log \cL_{2\a,\a})'(\lbd)\; =\; \frac{K_{\a-1}(2 \sqrt{\lbd})}{\sqrt{\lbd}K_{\a}(2 \sqrt{\lbd})}\; =\; \int_0^\infty \lpa\frac{1}{4\pi^2 u\,((J_\a(2\sqrt{u}))^2 + (Y_\a(2\sqrt{u}))^2)}\rpa \frac{du}{u+\lbd}
\end{equation}
where the second, non-trivial, equality follows from the main result of \cite{Gr} - see also \cite{I} for a simpler argument using the Perron-Stieltjes inversion formula and the Wronskian of Hankel functions. This shows that $\mu_{2\a,\a}$ has an explicit density $\varphi_{2\a, \a}$ which is expressed in terms of the classical Bessel functions $J_\a$ and $Y_\a:$
$$\varphi_{2\a, \a}(u)\; =\; \frac{1}{4\pi^2 u\,((J_\a(2\sqrt{u}))^2 + (Y_\a(2\sqrt{u}))^2)}\cdot$$

\begin{Rem} {\em For every $\a > 0, t\neq 0,$ the Laplace transform of $\GG_\a^{-\frac{1}{t}}$ is computed formally as
$$\esp[e^{-\lbd \GG_\a^{-\frac{1}{t}}}]\; =\; \frac{1}{t \Ga(\a)}\int_0^\infty x^{\a t -1} \exp\lpa - x^t -\frac{\lbd}{x}\rpa dx\; =\; \frac{Z^{\a t}_t(\lbd)}{t\Ga(\a)},$$
where $Z_\rho^\nu$ is the so-called Kr\"atzel function - see (1.7.42) in \cite{KST}. On the other hand, we know by Theorem 4 in \cite{BS2} and the discussion therebefore, that
$$\GG_\a^{-\frac{1}{t}}\; \in\;\cG\; \Leftrightarrow\;\GG_\a^{-\frac{1}{t}}\; \in\;\cI\; \Leftrightarrow\; t\ge -1.$$
This shows that $-(\log Z_\rho^\nu)'$ is the Stieltjes transform of a positive measure $\mu_{\rho,\nu}$ for all $\rho\ge -1$ and $\nu\rho \ge 0,$ and that it is not CM for $\rho < -1.$ The measure $\mu_{\rho,\nu}$ is not explicit in general, except for $\rho = 1$ by the preceding discussion and Formula (1.7.43) in \cite{KST}. The case $\rho =\nu > 0$ corresponds to the Fr\'echet random variable $\GG_1^{-\frac{1}{\rho}}$ and to our above special case $\a =1.$ It is also discussed in Section 3.4 of \cite{BS1} for $\rho =\nu \in (0,1),$ from the point of view of Bochner's subordination.}

\end{Rem}

\subsubsection{The case $m = 2$} As seen in the introduction, we have
$$\cL_{2, \a}(\lbd)\; =\; \frac{2 \nu^\nu}{\Ga(\nu)}\, \sqrt{\lbd} \, K_\nu (2\nu \lbd^{\frac{1}{2\nu}})$$
with the notation $\nu =\frac{1}{2-\a}\in (1/2, \infty).$ The same computation as above and the Perron-Stieltjes inversion formula lead to
$$-(\log \cL_{2,\a})'(\lbd)\; =\;\frac{\lbd^{\frac{1}{2\nu} -1} K_{\nu-1}(2\nu \lbd^{\frac{1}{2\nu}})}{K_\nu(2\nu \lbd^{\frac{1}{2\nu}})}\; =\; \int_0^\infty \lpa\frac{1}{2\pi\nu u}\, \Im \lpa \frac{z K_{\nu-1} (z)}{K_\nu (z)}\rpa\rpa \frac{du}{u+\lbd}$$
with $z = 2\nu (e^{\i\pi} u)^{\frac{1}{2\nu}}$, which shows a semi-explicit expression for the density $\varphi_{2, \a}$ of $\mu_{2,\a}.$

When $\nu > 1\Leftrightarrow \a \in (1,2),$ we have $\arg (z^2) \in (0,\pi)$ and we can apply the second equality in (\ref{PNT}) which holds on the complex plane cut along the negative real axis. This yields, after some algebra, an explicit integral representation connecting $\varphi_{2,\a}$ to $\varphi_{2\a, \a}:$
$$\varphi_{2, \a} (u) \; = \; \frac{1}{\nu^2} \int_0^\infty f_{\X_\nu}(\frac{u}{x}) \lpa x^{\frac{1}{\nu} -1} \varphi_{2\a, \a} (x^{\frac{1}{\nu}}) \rpa \frac{dx}{x}\,,$$
where $f_{\X_\nu}$ is the density of $\X_\nu\elaw \nu^{-2\nu}(\C_{\frac{1}{\nu}})^\nu$ and $\C_\mu$ is for every $\mu\in(0,1)$ the half-Cauchy random variable with density
$$\frac{\sin(\pi\mu)}{\pi\mu(x^2+2\cos(\pi\mu) x +1)}\cdot$$
Observe that since $\X_\nu\claw 1$ as $\nu\to 1,$ the above representation boils down to the tautological identity $\varphi_{2, 1} = \varphi_{2, 1}$ when $\nu = \a =1,$ a special case of Paragraph 3.2.1 above.

When $\nu \in(1/2, 1)\Leftrightarrow \a \in (0,1),$ we can write $z = \i Z$ with $Z = 2\nu e^{\i\pi(\frac{1}{2\nu}-\frac{1}{2})} u^{\frac{1}{2\nu}}$ such that  $\arg (Z^2) \in (0,\pi).$ By Formula 7.2.2(16) in \cite{EMOT}, we obtain
$$\varphi_{2, \a} (u) \; = \; \frac{-1}{2\pi\nu u}\, \Im \lpa \frac{Z H^{(2)}_{\nu-1} (Z)}{H^{(2)}_\nu (Z)}\rpa,$$
a complex expression which does not seem to lead to any particular real simplification.

\subsection{Some limit behaviours} In this last paragraph we briefly mention the limit behaviour of $\cG(m,\a)$ when the parameters $(m,\a)$ reach their admissibility boundary.

\begin{itemize}

\item When $\a\to 0 \Leftrightarrow a \to m,$ our main result shows immediately that
$$X_{m,\a}\; \claw\; 1,$$
an extension of the case $m = 1$ where it is obvious from (\ref{LaSta}) that $\Z_a\claw 1$ as $a\to 1.$

\item When $\a \to m\Leftrightarrow a \to 0,$ the second identity of our main result shows that
\begin{eqnarray*}
a X_{m, \a}^a\!\! &\!\!\elaw & \frac{\Gamma(a+1)}{\Gamma(m)} \times\lpa \frac{m}{a} \times \B_{1,\frac{m}{a}-1}\rpa^{-1}\times\lpa\prod_{n =0}^\infty \lpa\frac{m+n+1}{a+n+1}\rpa\B_{1+\frac{1+n}{a}, \frac{m}{a}-1} \rpa^{-1}\\
& \claw & \frac{1}{\Ga (m) \GG_1},
\end{eqnarray*}
where the convergence follows from (2.5) in \cite{LS2}. This is again an extension of the case $m=1$ where $\Z_a^{-a}\claw \GG_1$ as $a\to 0.$ See also \cite{C} for the asymptotic behaviour of real stable laws with small self-similarity parameter.

\end{itemize}

\begin{Rem}{\em When $m\to \infty$ and $\a$ or $a$ is fixed, putting together (\ref{MMS}) and (4.5) in \cite{BK} shows after some comparison with Theorem 1.4 and Remark 1.5 in \cite{NY} that $\cG(m,\a)$ exhibits a so-called mod-Gaussian convergence. We have not investigated the full details as yet, leaving them to further research.}
\end{Rem}

\bigskip
\noindent
{\bf Acknowledgement.}  The work of the first author is supported by King Saud University, Deanship of Scientific Research, College of Science Research Center.


\begin{thebibliography}{10}

\bibitem{Bar}
E.~W.~Barnes. The genesis of the double gamma function. {\em Proc. London Math. Soc.} {\bf 31}, 358-381, 1899.

\bibitem{BSW}
J.~Bernasconi, W.~R.~Schneider and W.~Wyss. Diffusion and hopping conductivity in disordered one-dimensional lattice systems. {\em Z. Phys. B} {\bf 37}, 175-184, 1980.

\bibitem{BY}
J.~Bertoin and M.~Yor. On the entire moments of self-similar Markov processes and exponential functionals. {\em Ann. Fac. Sci. Toulouse VI. S\'er. Math.} {\bf 11}, 33-45, 2002.

\bibitem{BK}
J.~Billingham and A.~C.~King. Uniform asymptotic expansions for the Barnes double gamma function. {\em Proc. Roy. Soc. London Ser. A} {\bf 453}, 1817-1829, 1997.

\bibitem{Bond}
L.~Bondesson. {\em Generalized Gamma convolutions and related classes of distributions and densities.} Lect. Notes Stat. {\bf 76}, Springer-Verlag, New York, 1992.

\bibitem{BJTP}
L.~Bondesson. A class of probability distributions that is closed with respect to addition as well as multiplication of independent random variables. {\em J. Theor. Probab.} {\bf 28} (3), 1063-1081, 2015.

\bibitem{BS1}
P.~Bosch and T.~Simon. On the self-decomposability of the Fr\'echet distribution.
{\em Indag. Math.} {\bf 24}, 626-636, 2013.

\bibitem{BS2}
P.~Bosch and T.~Simon. On the infinite divisibility of inverse Beta distributions. {\em Bernoulli} $\,$ {\bf 21} (4), 2552-2568, 2015.

\bibitem{BS3}
P.~Bosch and T.~Simon. A proof of Bondesson's conjecture on stable densities. {\em Ark. Mat.} {\bf 54}, 31-38, 2016.

\bibitem{Br}
B.~L.~J.~Braaksma. Asymptotic expansions and analytic continuations for a class of Barnes-integrals. {\em Compos. Math.} {\bf 15}, 239-341, 1963.

\bibitem{C}
N.~Cressie. A note on the behaviour of the stable distributions for small index $\a.$ {\em Z. Wahrsch. verw. Gebiete} {\bf 33}, 61-64, 1975.

\bibitem{CSY}
M.~Cs\"org\"o, Z.~Shi and M.~Yor. Some asymptotic properties of the local
time of the uniform empirical process. {\em Bernoulli} {\bf 5} (6), 1035-1058, 1999.

\bibitem{ET}
A.~Erd\'elyi and F.~G.~Tricomi. The asymptotic expansion of a ratio of Gamma functions. {\em Pacific J. Math.} {\bf 1} (1), 133-142, 1951.

\bibitem{EMOT}
A.~Erd\'elyi, W. Magnus, F.~Oberhettinger and F.~G.~Tricomi. {\em Higher transcendental functions. Vol. II.} McGraw-Hill, New-York, 1953.

\bibitem{FGD}
P.~Flajolet, X.~Gourdon and P.~Dumas. Mellin transforms and asymptotics: Harmonic sums. {\em Theoret. Comput. Sci.} {\bf 144}, 3-58, 1995.

\bibitem{Gr}
E.~Grosswald. The Student $t-$distribution of any degree of freedom is infinitely divisible. {\em Z. Wahrsch. verw. Gebiete} {\bf 36}, 103-109, 1976.

\bibitem{IL}
I.~A.~Ibragimov and Yu.~V.~Linnik. {\em Independent and stationary sequences of random variables.} Wolters-Noordhoff, Groningen, 1971.

\bibitem{I}
M.~H.~E.~Ismail. Bessel functions and the infinite divisibility of the Student $t-$distribution. {\em Ann. Probab.} {\bf 5}, 582-585, 1977.

\bibitem{JRY}
L.~F.~James, B.~Roynette and M.~Yor. Generalized gamma convolutions, Dirichlet means, Thorin measures, with explicit examples. {\em Probab. Surv.} {\bf 5}, 346-415, 2008.

\bibitem{J}
S.~Janson. Moments of Gamma type and the Brownian supremum process area. {\em Probab. Surveys} {\bf 7}, 1-52, 2010.

\bibitem{jedidi}  
W.~Jedidi.  Stable processes, mixing, and distributional properties, Part II. {\em
Teor. Veroyatn. Primen.} {\bf  53} (1), 81-105, 2009.

\bibitem{KST}
A.~A.~Kilbas, H.~M.~Srivastava and J.~J.~Trujillo. {\em Theory and applications of fractional differential equations.} North-Holland, Amsterdam, 2006.

\bibitem{KP}
A.~Kuznetsov and J.-C.~Pardo. Fluctuations of stable processes and exponential functionals of hypergeometric L\'evy processes. {\em Acta Appl. Math.} {\bf 123}, 113-139, 2013.

\bibitem{LS2}
J.~Letemplier and T.~Simon. On the law of homogeneous stable functionals. Preprint. {\tt arXiv:1510.01625}

\bibitem{NY}
A.~Nikeghbali and M.~Yor. The Barnes $G-$function and its relations with sums and products of generalized Gamma convolution variables. {\em Elec. Comm. in Probab.} {\bf 14}, 396-411, 2009.

\bibitem{Pk97}
A.~G.~Pakes. Characterization by invariance under length-biasing and random scaling. {\em J. Statist. Plann. Inference} {\bf 63}, 285-310, 1997.

\bibitem{Pk14}
A.~G.~Pakes. On generalized stable and related laws. {\em J. Math. Anal. Appl.} {\bf 411}, 201-222, 2014.

\bibitem{PS3}
P.~Patie and M.~Savov. Spectral expansions of non-self-adjoint generalized Laguerre semigroups. Preprint. {\tt arXiv:1506.01625}

\bibitem{WRS}
W.~R.~Schneider. Generalized one-sided stable distributions. In: {\em Stochastic processes - mathematics and physics II}, Lect. Notes Math. {\bf 1250}, 269-287, 1987.

\bibitem{SVH}
F.~W.~Steutel and K.~van Harn. {\em Infinite divisibility of probability distributions on the real line.} Marcel Dekker, New-York, 2003.

\bibitem{Wb}
R.~Webster. Log-convex solutions to the functional equation $f(x+1) = g(x)f(x): \Gamma-$type functions. {\em J. Math. Anal. Appl.} {\bf 209}, 605-623, 1997.

\end{thebibliography}
\end{document}